

\magnification\magstep1
\vsize=24.0truecm
\baselineskip13pt
\def\knep{\noalign{\vskip2pt}}

\def\eps{\varepsilon}
\def\half{\hbox{$1\over2$}}

\font\csc=cmcsc10
\font\cyr=wncyr10
\font\bigbf=cmbx12
\font\bigcyr=wncyr10 at 12pt
\font\bigrm=cmr10 at 12pt

\def\j3ord{\rm\u{\cyr i}}
\def\j3bold{\bf\u{\bigcyr i}}

\def\cstok#1{\leavevmode\thinspace\hbox{\vrule\vtop{\vbox{\hrule\kern1pt
        \hbox{\vphantom{\tt/}\thinspace{\tt#1}\thinspace}}
        \kern1pt\hrule}\vrule}\thinspace} 
\def\firkant{\cstok{\phantom{$\cdot$}}} 

\centerline{\bigbf On the time a diffusion process spends along a line}

\medskip 
\centerline{\bigrm Nils Lid Hjort and 
	{\bigcyr Khasp1minski{\j3bold} Rafail Zalmanovich}}

\medskip
\centerline{\bigrm University of Oslo and}
\centerline{\bigcyr Institut Problem Peredachi Informatsii, Moskva}

\smallskip
\centerline{\tt -- October 1992 --}

\smallskip
{{\smallskip\narrower\baselineskip12pt\noindent
{\csc Abstract.} 
For an arbitrary diffusion process $X$ with 
time-homogeneous drift and variance parameters 
$\mu(x)$ and $\sigma^2(x)$, 
let $V_\eps$ be $1/\eps$ times the total time
$X(t)$ spends in the strip 
$[a+bt-\half\eps,a+bt+\half\eps]$. 
The limit $V$ as $\eps\rightarrow0$ 
is the full {half}line version of the 
local time of $X(t)-a-bt$ at zero, 
and can be thought of as 
the time $X$ spends along the straight line $x=a+bt$.  
We prove that $V$ is either infinite with probability 1 or 
distributed as a mixture of an exponential 
and a unit point mass at zero, 
and we give formulae for the parameters of this distribution 
in terms of 
$\mu(.)$, $\sigma(.)$, $a$, $b$, and the starting point $X(0)$.   
The special case of a Brownian motion 
is studied in more detail, 
leading in particular to 
a full process $V(b)$ with continuous sample paths
and exponentially distributed marginals. 
This construction leads 
to new families of bivariate and multivariate exponential distributions.
Truncated versions of such `total relative time' variables are 
also studied. A relation is pointed out 
to a second order asymptotics problem in statistical estimation theory,
recently investigated in Hjort and Fenstad (1992a, 1992b). 

\smallskip\noindent
{\csc Key Words:} {\sl Brownian motion, 
diffusion process, 
exponential process, 
local time, 
multivariate exponential distribution, 
second order asymptotics for estimators }
\smallskip}}

\medskip 
{\bf 1. Introduction and summary.}
Consider a time-homogeneous diffusion process $X$ with 
$dX(t)=\mu(X(t))\,dt+\sigma(X(t))\,dW(t)$, 
using $W(.)$ to denote a standard Brownian motion. 
In other words, $X$ is a Markov process with continuous paths 
and with the property that $X(t+h)-X(t)$, for given $X(t)=x$,
has mean value $\mu(x)h+o(h)$ and variance $\sigma^2(x)h+o(h)$.
Consider 
$$V_\eps={1\over \eps}\int_0^\infty
	 I\{|X(t)-a-bt|\le\half\eps\}\,dt, \eqno(1.1) $$
the total amount of time spent by the process in 
the narrow strip $[a+bt-\half\eps,a+bt+\half\eps]$, 
divided by $\eps$. 
We show in Section 2 that this variable has a well-defined
limit $V$, and find its distribution, in terms of 
$a$, $b$, $\mu(.)$, $\sigma(.)$, and the starting point $X(0)=x$.
Under certain conditions the variable is infinite almost surely,
and in the opposite case the variable 
is distributed as a mixture of an exponential and a unit point mass
at zero. Explicit formulae for the parameters of this
distribution are also found.
The variable is a pure exponential in the case $X(0)=a$.   
The simplicity of the result is remarkable,
in view of the large class of diffusion processes; 
in particular $X$ can have both Gau\ss ian and non-Gau\ss ian sample paths.  

The $V$ variable can be thought of as the total
relative time the $X(.)$ process spends along the line $x=a+bt$,
and is related to what is sometimes called the local time at zero
of the $X(t)-a-bt$ process. 
Usually such local times are studied and used 
for a limited time interval $[0,\tau]$ only, however. 

A special case of the construction above is that of $a=0$, giving 
$$V_\eps(b)={1\over \eps}\,{\rm measure}\{t\ge0\colon W(t)
	\in [bt-\half\eps,bt+\half\eps]\}. \eqno(1.2)$$
That the limit $V(b)$, the time Brownian motion spends along $x=bt$, 
is simply exponential with parameter $|b|$, 
follows from the general result of Section 2, 
but is proved in a more direct fashion in Section 3, 
using moment convergence. 
This second approach lends itself more easily
to the simultaneous study of several relative times.
In Section 4 we prove full process convergence 
of $\{V_\eps(b)\colon b\not=0\}$ towards 
a $\{V(b)\colon b\not=0\}$ with continuous sample paths and  
with exponentially distributed marginals. 
Its covariance and correlation structure is also found. 
In particular this construction leads to new families of 
bivariate and multivariate exponential distributions.

In Section 5 a more general variable $V(c,b)$ is studied,
defined as the total relative time during $[c,\infty)$ that  
$W(t)$ spends along $bt$. The distribution of $V(c,b)$ is again 
a mixture of an exponential and a unit point mass at zero. 
A simple consequence
of this result is a rederivation of a well known formula for
the distribution of the maximum of Brownian motion over an interval. 
Finally some supplementing results and remarks are given in Section 6.
In particular some consequences for empirical partial sum processes
are briefly discussed. 

A certain second order asymptotics
problem in statistical estimation theory led by serendipity to 
the present study on total relative time variables for Brownian motion.  
Suppose $\{\theta_n\colon n\ge1\}$ is an estimator sequence
for a parameter $\theta$, where $\theta_n$ is based on the first 
$n$ data points in an i.i.d.~sequence, and consider 
$Q_\delta$, the number of times, among $n\ge c/\delta^2$, 
where $|\theta_n-\theta|\ge\delta$. 
Almost sure convergence (or strong consistency) of $\theta_n$ 
is equivalent to saying that $Q_\delta$ is almost surely 
finite for every $\delta$, and it is natural to inquire 
about its size. A particular result of Hjort and Fenstad (1992a, Section 7)
is that under natural conditions, which include the existence of 
a normal $(0,\sigma^2)$ limit for $\sqrt{n}(\theta_n-\theta)$,
$$\delta^2Q_\delta\rightarrow_d 
	Q=Q(c,1/\sigma)
	 =\int_c^\infty I\{|W(t)|\ge t/\sigma\}\,dt \eqno(1.3)$$
as $\delta\rightarrow0$. 
If $\{\theta_{n,1}\}$ and $\{\theta_{n,2}\}$ are 
first order equivalent estimator sequences,  
with the same $N(0,\sigma^2)$ limit for 
$\sqrt{n}(\theta_{n,j}-\theta)$, 
and $Q_{\delta,j}$ is the number of $\delta$-misses for sequence $j$, 
then $Q_{\delta,1}/Q_{\delta,2}\rightarrow1$ and 
$\delta^2(Q_{\delta,1}-Q_{\delta,2})\rightarrow0$ in probability. 
One way of distinguishing between the two estimation methods 
is by studying second order aspects of $Q_{\delta,1}-Q_{\delta,2}$. 
It turns out that $\delta$ times this difference in
typical cases has a limit distribution which is that of a constant times 
$V(c,1/\sigma)-V(c,-1/\sigma)$, or times the simpler 
$V(1/\sigma)-V(-1/\sigma)$ if $c=c(\delta)$ is allowed to decrease to zero
in the definition of $Q_{\delta,j}$. 
Note the connection from $Q(c,1/\sigma)$ of (1.3) 
to $V(c,\pm 1/\sigma)$.  
Some further details are in 6C in the present paper,
while further background and discussion can be found 
in Hjort and Fenstad (1992a, 1992b). 

\bigskip
{\bf 2. The time X spends along a straight line.} 
In 2A we solve the problem for the time spent along
a line parallel to the time axis. 
This rather immediately leads to the
more general solution, which is presented in 2B.

{\sl 2A. The time X spends along a horizontal line.}
Let $X(t)$ be as in the introductory paragraph,
with continuous and positive diffusion function $\sigma(x)$ 
and continuous drift function $\mu(x)$. For a temporarily 
fixed $a$, define 
$$s(y)=\exp\Bigl\{-\int_a^y {2\mu(x)\over \sigma^2(x)}\Bigr\}, $$
also for negative $y$. The function 
$S(z)=\int_a^{z} s(y)\,dy$, or any linear translation thereof, 
is often called the scale function of the diffusion process. 
Two important quantities are 
$$k_+(a)=\int_a^\infty s(y)\,dy 
	\quad {\rm and} \quad 
  k_-(a)=\int_{-\infty}^a s(y)\,dy. \eqno(2.1)$$
It is known that if $k_+(a)$ is finite, then there is 
a positive probability for the process to drift off towards $+\infty$, 
and vice versa; and similarly the finiteness of $k_-(a)$ 
corresponds exactly to there being a positive probability for
drifting off towards $-\infty$. 
See for example Karlin and Taylor (1981, Chapter 15.6). 
If in particular both integrals are infinite then 
the process is recurrent and visits the line $x=a$ 
an infinite number of times. 

The current object of interest is 
$$V_\eps={1\over \eps}\int_0^\infty 
	I\{|X(t)-a|\le\half\eps\}\,dt. \eqno(2.2)$$ 
Let $V_\eps(\tau)$ be defined similarly, but for 
the interval $[0,\tau]$ only. 
This is the so-called local
time at zero process, Paul L\'evy's `mesure de voisinage',
for $X(t)-a$; see for example Karlin and Taylor 
(1981, Chapter 15.12) and It\^o and McKean (1979, Chapters 2 and 6).
It is a `remarkable and recondite fact', 
to quote Karlin and Taylor,
that the limit $V(\tau)$ of $V_\eps(\tau)$ as $\eps\rightarrow0$
exists for almost every sample path 
[that is, $V_\eps(\tau,\omega)$ converges to a well-defined
$V(\tau,\omega)$ for each $\omega$ in a subset of probability 1 
of the underlying probability space]. 
It follows from this local time theory that 
$V_\eps=V_\eps(\infty)$ converges to a well-defined 
$V=V(\infty)$ too, with probability 1. 
We think of $V$ as the total relative time $X$ spends along the
line $x=a$. 

In the following we are able to find the exact distribution
of $V$. The arguments we shall use actually show 
convergence in distribution of $V_\eps$ to $V$ directly,
that is, we do not need or use the somewhat 
sophisticated local time theory or the almost sure 
pathwise existence of $V$ to prove that $V_\eps$ 
has the indicated limit distribution as $\eps$ goes to zero.  

If $\alpha$ is positive, 
write $V\sim{\rm Exp}(\alpha)$ for the exponential distribution with 
density $g(v)=\alpha e^{-\alpha v}$ for $v\ge0$. 
It has mean $1/\alpha$ and Laplace transform 
$E\exp(-\lambda V)=\alpha/(\alpha+\lambda)$. 

{{\smallskip\sl
{\csc Theorem 1.} 
Assume that the $X$ process starts at $X(0)=a$. 
If $k_+(a)$ and $k_-(a)$ are both infinite, then $V=\infty$ 
with probability one. 
Otherwise the limit $V$ of $V_\eps$ 
is exponentially distributed with parameter 
$\alpha(a)=\half\sigma^2(a)\{1/k_+(a)+1/k_-(a)\}$. 	 
\smallskip}}

{\csc Proof:}
That $V_\eps$ goes almost surely to infinity 
when both integrals are infinite follows from the 
theory of Karlin and Taylor (1981, Chapter 15.6).
This is connected to the recurrency phenomenon 
mentioned after (2.1) above. 

Suppose next that both $k_+(a)$ and $k_-(a)$ are finite. 
For a fixed positive $\lambda$, study the function 
$$u_{\lambda,\eps}(x)
	=E_x\exp(-\lambda V_\eps)	
	=E_x\exp\bigl\{-\lambda\int_0^\infty  f_\eps(X(t))\,dt\bigr\}, $$
where the subscript $x$ here and below 
means that the expectation is conditional on starting point $X(0)=x$, 
and where $f_\eps(x)=\eps^{-1}I\{|x-a|\le\half\eps\}$. 
General results for 
diffusion processes imply that 
the $u_{\lambda,\eps}$ function has two piecewise 
continuous derivatives and satisfies   
$$\half\sigma^2(x)u_{\lambda,\eps}''(x)
	+\mu(x)u_{\lambda,\eps}'(x)
	-\lambda f_\eps(x)u_{\lambda,\eps}(x)=0, $$
see for example the theory developed by 
Karlin and Taylor (1981, Chapter 15.3). 
Integrating from $a-\half\eps$ to $a+\half\eps$ 
and letting $\eps\rightarrow0$
shows that a solution $u_\lambda(x)$ to the limit problem
must satisfy 
$$u_\lambda'(a+)-u_\lambda'(a-)
	=2\lambda u_\lambda(a)/\sigma^2(a). \eqno(2.3)$$


Now let $w(x,a)$ be the probability that the process after start in 
$X(0)=x$ succeeds in reaching the level $x=a$ in a finite amount
of time. If this happens then $V$ starting from $x$ is equal 
in distribution to a $V$ starting from $a$, 
because of the strong Markov property and the postulated time-homogeneity. 
And if it does not happen then $V=0$. Hence 
$$u_\lambda(x)=E_xe^{-\lambda V}
	=w(x,a)\,E_ae^{-\lambda V}+\{1-w(x,a)\}\,Ee^{-0}
	=w(x,a)u_\lambda(a)+1-w(x,a). \eqno(2.4) $$ 
This equation is also reached if one more carefully 
starts with $V_\eps$-equations and then lets $\eps\rightarrow0$. 
But $w(x,a)$ can be found explicitly, since it satisfies
$\half\sigma^2(x)w''(x,a)+\mu(x)w'(x,a)=0$ with 
boundary conditions $w(-\infty,a)=0$, $w(a,a)=1$, $w(\infty,a)=0$.
Differentiation here is w.r.t.~$x$ and $a$ is still fixed.  
The solution is 
$$w(x,a)=\cases{
k_+(x)/k_+(a) 
	&if $x\ge a$, \cr
	\knep
k_-(x)/k_-(a) 
	&if $x\le a$, \cr} \eqno(2.5)$$
in terms of the (2.1) functions. 	
In particular, 
$w'(a+,a)=-1/k_+(a)$ and $w'(a-,a)=1/k_-(a)$ in terms of the 
transience determining quantities (2.1).  
This can now be used in (2.3) to make (2.4) more explicit:
$$u_\lambda'(a+)=w'(a+,a)\{u_\lambda(a)-1\}
	\quad {\rm and} \quad 
  u_\lambda'(a-)=w'(a-,a)\{u_\lambda(a)-1\} $$
lead to $\{-1/k_+(a)-1/k_-(a)\}\{u_\lambda(a)-1\}
	=2\lambda u_\lambda(a)/\sigma^2(a)$.   
And solving this produces in the end 
$$u_\lambda(a)={1/k_+(a) +1/k_-(a)
	\over 1/k_+(a) +1/k_-(a) + 2\lambda/\sigma^2(a)}
	={\alpha(a)\over \alpha(a)+\lambda}, $$
with the $\alpha(a)$ parameter as given in the theorem. 

Assume next that $k_+(a)$ is finite but $k_-(a)$ infinite.
This case can be handled very much as the previous one. 
Now $+\infty$ is attracting but $-\infty$ is not, and
the boundary conditions for $w(x,a)$ become 
$w(-\infty,a)=1$, $w(a,a)=1$, $w(\infty,a)=0$, giving 
as solution  
$$w(x,a)=\cases{
k_+(x)/k_+(a) 
	&if $x\ge a$, \cr
1	&if $x\le a$. \cr} \eqno(2.6)$$
(2.3) and (2.4) are still valid, and we find 
after similar arguments that $V$ is exponential with 
parameter $\alpha(a)=\half\sigma^2(a)/k_+(a)$. 
The final case of $k_+(a)$ infinite and $k_-(a)$ finite 
is handled simiarly. 
\firkant 
 
\smallskip
The proof actually gives the distribution of $V$ for
an arbitrary starting point $x$, namely 
$$V|\{X(0)=x\}\sim w(x,a)\,{\rm Exp}(\alpha(a))
	+\{1-w(x,a)\}\,\delta_0, \eqno(2.7)$$
in which $\delta_0$ is a unit point mass at zero. 
The weight $w(x,a)$ here has a direct probabilistical 
interpretation, and is given in (2.5) for the case of 
two attracting boundaries and in (2.6) for the case 
of only $+\infty$ attracting, with a similar
modification for the case of $k_+(a)$ infinite but
$k_-(a)$ finite. 
In the case of (2.6) 
we see that $V\sim{\rm Exp}(\alpha(a))$ for any 
starting point to the left of $a$.

\smallskip
{\sl 2B. The time X spends along a general line from a general
starting point.} 
The generalisation to a result about the (1.1) variable 
is now immediate. 
Just consider the new process $X^*(t)=X(t)-bt$, 
which is a diffusion with $\mu^*(x)=\mu(x)-b$ 
and the same $\sigma(x)^2$. The previous result is
valid for the time $X^*(t)$ spends along the horizontal 
line $x^*=a$. We need 
$s_{a,b}(y)=\exp[-\int_a^y2\{\mu(x)-b\}/\sigma^2(x)\,dx]$ as well as 
$$k_+(a,b)=\int_a^\infty s_{a,b}(y)\,dy
	\quad {\rm and} \quad 
  k_-(a,b)=\int_{-\infty}^a s_{a,b}(y)\,dy. \eqno(2.8)$$
We find the following. 

{{\smallskip\sl
{\csc Theorem 2.} 
Let the process start at $X(0)=x$, and suppose 
one or both of the two integrals (2.8) are finite. 
Then $V_\eps$ of (1.1) converges in distribution to 
the mixture 
$$V|\{X(0)=x\}\sim w(x,a,b)\,{\rm Exp}(\alpha(a,b))
	+\{1-w(x,a,b)\}\,\delta_0 \eqno(2.9) $$
of an exponential and a unit point mass at zero. Here 
$$\alpha(a,b)=\half\sigma^2(a)\{k_+(a,b)^{-1}+k_-(a,b)^{-1}\} $$
and 
$$w(x,a,b)=\cases{
\int_x^\infty s_{a,b}(y)\,dy\big/k_+(a,b) 
	&if $x\ge a$, \cr
	\knep
\int_{-\infty}^{x} s_{a,b}(y)\,dy\big/k_-(a,b) 
	&if $x\le a$, \cr} $$
if both denominators are finite. If one of them is infinite,
replace the corresponding ratio with 1. 
\smallskip}}
 
{\sl 2C. Example: Total time for Brownian motion.}
Let us apply the general theorem to the case of 
$X=W$, the standard Brownian motion process, which has
$\mu(x)=0$ and $\sigma(x)=1$. We allow an arbitrary starting point
$W(0)=x$. Take $b$ positive and 
consider the total relative time $V_\eps$ of (1.1). Then
$$V_\eps|\{W(0)=x\}\rightarrow V\sim\cases{
{\rm Exp}(b) &if $x\ge a$, \cr
e^{-2b(a-x)}\,{\rm Exp}(b)+\{1-e^{-2b(a-x)}\}\,\delta_0
	     &if $x\le a$. \cr} \eqno(2.10)$$
There is a symmetric result for negative $b$, involving
an exponential with parameter $|b|$. Notice in particular 
that $V\sim{\rm Exp}(|b|)$ when the starting point is $a$. 
And when $b=0$ then $V$ is infinite with probability one;
see 6A for a more informative result. 

%
%
%
%
%

\bigskip
{\bf 3. Moment convergence proof.}
In the following we stick to the Brownian motion,
and for simplicity take it to start at $W(0)=0$. 
For $b\not=0$, let us consider 
$$V_\eps(b)={1\over \eps}\int_0^\infty I\{W(t)\in bt\pm\half\eps\}\,dt $$
of (1.2) in more detail. That 
$$V_\eps(b)\rightarrow_d V(b)\sim {\rm Exp}(|b|) \eqno(3.1)$$
is already a consequence of the general theorem,
and indeed a special case of (2.10) above.  
We now offer a different proof,  
by demonstrating appropriate convergence of all moments. 
This is sufficient since the exponential distribution is 
determined by its moment sequence. 
In addition to having some independent merit this proof
lends itself more easily to the study of simultaneous convergence
aspects; see Section 4. 

For the first moment, observe that 
$$EV_\eps(b)={1\over\eps}\int_0^\infty 
	{\rm Pr}\{bt-\half\eps\le W(t)\le bt+\half\eps\}\,dt
	=\int_0^\infty f_t(bt)\,dt+O(\eps), \eqno(3.2)$$
where $f_t(x)=\phi(x/\sqrt{t})/\sqrt{t}$ is the density function 
for $W(t)$. Accordingly $EV_\eps(b)$ goes to 
$\int_0^\infty \phi(b\sqrt{t})/\sqrt{t}\,dt=1/|b|$.
Next consider the $p$-th moment. One finds 
$$\eqalign{
EV_\eps(b)^p&={1\over \eps^p}\int_0^\infty\cdots\int_0^\infty
	{\rm Pr}\{W(t_1)\in bt_1\pm\half\eps,\ldots,
		  W(t_p)\in bt_p\pm\half\eps\}\,dt_1\cdots dt_p \cr
	&=p!\int\cdots\int_{t_1<\cdots<t_p} 
	 f_{t_1,\ldots,t_p}(bt_1,\ldots,bt_p)\,dt_1\cdots dt_p+O(\eps), \cr}$$
where $f_{t_1,\ldots,t_p}(x_1,\ldots,x_p)$ is 
the density function of $(W(t_1),\ldots,W(t_p))$.
By the Gau\ss ian and Markovian 
properties of $W(.)$ this density can in fact be written
$$\phi\Bigl({x_1\over\sqrt{t_1}}\Bigr){1\over\sqrt{t_1}}
	 \phi\Bigl({x_2-x_1\over\sqrt{t_2-t_1}}\Bigr){1\over\sqrt{t_2-t_1}}
	 \cdots 
	 \phi\Bigl({x_p-x_{p-1}\over\sqrt{t_p-t_{p-1}}}\Bigr)
       				{1\over\sqrt{t_p-t_{p-1}}} \eqno(3.3)$$
when $t_1<\cdots <t_p$. To carry out the $p$-dimensional integration,
insert $(bt_1,\ldots,bt_p)$ for $(x_1,\ldots,x_p)$, and 
transform to new variables $u_1=t_1$, $u_i=t_i-t_{i-1}$ for $i=2,\ldots,p$.
The result is then that 
$$EV_\eps(b)^p\rightarrow 
	p!\int_0^\infty\cdots\int_0^\infty 
	{\phi(b\sqrt{u_1})\over \sqrt{u_1}}\cdots
	{\phi(b\sqrt{u_p})\over \sqrt{u_p}}\,du_1\cdots du_p
	=p!\,(1/|b|)^p$$
for each $p$. But this is manifestly 
the moment sequence of ${\rm Exp}(|b|)$, proving (3.1). \firkant

The case of $b=0$ is different, since $W$ spends an infinite amount 
of time along the time axis. An interesting $|N(0,1)|$ 
limit result for the relative time in $\pm\eps$ during $[0,T]$ 
is in 6A below. 

\bigskip
{\bf 4. The exponential process.}
We have seen that $V_\eps(b)$ goes to an exponentially distributed
$V(b)$, and in the same manner we should find bivariate and multivariate
exponential distributions by considering 
two or more $b$'s at the same time. This requires verification of 
simultaneous convergence in distribution of $(V_\eps(b),V_\eps(c))$ 
and similar quantities. This section indeed demonstrates 
process convergence of $V_\eps(.)$ to $V(.)$, 
and studies some of the properties of the limiting process.

\smallskip
{\sl 4A. Process convergence.} The first main result is as follows.

\smallskip
{\csc Theorem.} 
{\sl There is a well-defined stochastic process 
$V=\{V(b)\colon b\not=0\}$ 
with exponentially distributed marginals 
and with the property that 
$(V(b_1),\ldots,V(b_n))$
is the limit in distribution of 
$(V_\eps(b_1),\ldots,V_\eps(b_n))$ 
for each finite set of non-null indexes $b_i$. 
There exists a version of \/$V$\/with continuous paths, and  
$V_\eps(.)$ converges to $V(.)$ in the uniform topology 
on the $C$-space $C[b_0,b_1]$ of continuous functions on $[b_0,b_1]$, 
for each interval not containing zero.}

{\csc Proof:} 
Consider two rays $bt$ and $ct$ and their associated 
total relative time variables $(V_\eps(b),V_\eps(c))$.
Using the Cram\'er--Wold theorem in conjunction 
with the moment convergence method we see that convergence of 
$EV_\eps(b)^pV_\eps(c)^q$ to the appropriate limit, for each $p$ and $q$, 
is sufficient.
But this can be proved by slight elaborations on 
the techniques of Section 3. By Fubini's theorem 
$$\eqalign{
V_\eps(b)^pV_\eps(c)^q
	 ={1\over \eps^{p+q}}
	 \int_0^\infty\cdots\int_0^\infty
	 I\{&W(s_1)\in bs_1\pm\half\eps,\ldots,W(s_p)\in bs_p\pm\half\eps, \cr
	    &W(t_1)\in ct_1\pm\half\eps,\ldots,W(t_q)\in ct_q\pm\half\eps\}
	  \,ds_1\cdots dt_q, \cr}$$
and its expected value is seen to converge to 
$$EV(b)^pV(c)^q
	=\int_0^\infty\cdots\int_0^\infty
	 f_{s_1,\ldots,s_p,t_1,\ldots,t_q}
	  (bs_1,\ldots,bs_p,ct_1,\ldots,ct_q) 
	  \,ds_1\cdots dt_q \eqno(4.1)$$
by Lebesgue's theorem on dominated convergence. 
Note that the integral is over all of $[0,\infty)^{p+q}$ 
and that a simple expression like (3.3) for the density 
of $(W(s_1),\ldots,W(t_q))$ is only valid when the
time-points are ordered, so the factual integration in 
(4.1) is difficult to carry through 
(but possible; see 4C below).  
What is important at the moment
is however the mere existence of this and all other similar limits of 
product moments for the $V_\eps(.)$-process. 
We may conclude that all finite-dimensional distributions 
converge to well-defined limits. 
That these finite-dimensional distributions also 
constitute a Kolmogorov-consistent system is a by-product 
of the tightness condition verified below. 

The $V_\eps(.)$-process has continuous paths in $b\not=0$ for each $\eps$, 
since $W(.)$ is continuous. 
In order to prove process convergence on $C[b_0,b_1]$ for
a given interval we need to demonstrate tightness 
of the $\{V_\eps(.)\}$ family as $\eps$ goes to zero. 
Note that if $V_\eps^*(b)=V_\eps(-b)$, then
the processes $V_\eps^*(.)$ and $V_\eps(.)$ 
have identical distributional characteristics, 
so it suffices to consider the positive part of the process. 
Following results in Shorack and Wellner (1986, page 52) 
it is enough to verify that 
$$\limsup_{\eps\rightarrow0}E\{V_\eps(b+h)-V_\eps(b)\}^4\le Kh^2
	\quad {\rm for\ some\ }K, \eqno(4.2)$$
for all $h\ge0$ and for all $b$ with $b$ and $b+h$ in $[b_0,b_1]$,
where $0<b_0<b_1$.  
By the arguments for finite-dimensional convergence used above 
the left hand side of (4.2) is equal to 
$m_b(h)=E\{V(b+h)-V(b)\}^4$. This is seen to be a 
smooth function of $h$ with finite derivatives at zero. 
Ingenious and rather elaborate Taylor expansion arguments 
can in fact be furnished to prove that 
$$\eqalign{
EV(b+h)^4      &=(24/b^4)\{1-4\delta+O(\delta^2)\}, \cr
EV(b)V(b+h)^3  &=(24/b^4)\{1-6\delta+O(\delta^2)\}, \cr
EV(b)^2V(b+h)^2&=(24/b^4)\{1-6\delta+O(\delta^2)\}, \cr
EV(b)^3V(b+h)  &=(24/b^4)\{1-4\delta+O(\delta^2)\}, \cr
EV(b)^4        &=24/b^4, 			    \cr}$$
where $\delta=h/b$, so that 
$m_b(h)=K_2(b)h^2+K_3(b)h^3+\ldots$,
for local constants $K_j(b)$ that are continuous functions of $b$
(as long as $b\not=0$).   
This is dominated by a common $Kh^2$ for all $b$ and $b+h$ 
in the interval under consideration. This verifies (4.2), 
and incidentally at the same time 
verifies the so-called Kolmogorov condition
for almost sure continuity of the sample paths, see Shorack and Wellner
(1986, Chapter 2, Section 3). 

Using the moment formula in 4C below one may in fact 
calulate the left hand side of (4.2) explicitly, 
and a fair amount of analysis leads to 
$m_b(h)=24\cdot352\,h^2/b^6+O(h^3)$.
The proof above circumvented the need for 
information on this level of detail, however. \firkant

\smallskip 
{\sl 4B. Dependence structure.}
In order to investigate this to some extent we 
calculate covariances and correlations. 
Let $0<b<c$ and $-c<0<d$. Then 
$${\rm cov}\{V(b),V(c)\}={1\over c}{1\over 2c-b}
	\quad {\rm and} \quad 
  {\rm cov}\{V(-c),V(d)\}
	={1\over d}{1\over c+2d}
		+{1\over c}{1\over 2c+d}-{1\over cd}.\eqno(4.3)$$
To prove this, consider the case of two positive parameters. 
Then by previous arguments 
$$\eqalign{
EV_\eps(b)V_\eps(c)
&=\int_0^\infty\int_0^\infty
	{\rm Pr}\{W(s)\in bs\pm\half\eps,\,
	W(t)\in ct\pm\half\eps\}\,dsdt/\eps^2 \cr
&\rightarrow\int\int_{s<t}\bigl[f_{s,t}(bs,ct)
		+f_{s,t}(cs,bt)\bigr]\,dsdt \cr
&=\int\int_{s<t}\Bigl[\phi(b\sqrt{s})\phi\Bigl({ct-bs\over \sqrt{t-s}}\Bigr)
	+\phi(c\sqrt{s})\phi\Bigl({bt-cs\over \sqrt{t-s}}\Bigr)\Bigr]
	{1\over \sqrt{s}}{1\over \sqrt{t-s}}\,dsdt, \cr}$$
where (3.3) is used again. Now transform first to $(s,u)=(s,t-s)$ 
and then to $(x,y)=(\sqrt{s},\sqrt{u})$, to get
$$4\int_0^\infty\int_0^\infty
  \Bigl[\phi(bx)\phi\Bigl({cy^2+(c-b)x^2\over y}\Bigr)
	+\phi(cx)\phi\Bigl({by^2-(c-b)x^2\over y}\Bigr)\Bigr]\,dxdy.$$
The rest of the calculation is carried out using the formula  
$\int_0^\infty\exp\{-\half(k^2y^2+l^2/y^2)\}\,dy
	=\half\sqrt{2\pi}(1/k)\exp(-kl)$. 
This formula can be proved by clever but elementary integrations, 
and is valid for positive $k$ and $l$. One finds
$${4\over 2\pi}{\sqrt{2\pi}\over2}
	\int_0^\infty\Bigl[{1\over c}\exp\{-\half(b^2+4c(c-b))x^2\}
	+{1\over b}\exp(-\half c^2x^2)\Bigr]\,dx
	={1\over c}{1\over 2c-b}+{1\over b}{1\over c}.$$
The first formula in (4.3) follows from this,
and the other case is handled similarly. \firkant

\smallskip
It is convenient to give formulae (4.3) in another form,
using $(b,c)=(b,b+h)$ in the first case and $(-c,d)=(-c,kc)$ in the second.
Then 
$${\rm cov}\{V(b),V(b+h)\}={1\over b+h}{1\over b+2h}
	\quad {\rm and} \quad
  	{\rm cov}\{V(-c),V(kc)\}={1\over c^2}{-3\over (k+2)(2k+1)},$$
and the correlation coefficients become
$${\rm corr}\{V(b),V(b+h)\}={b\over b+2h}
	\quad {\rm and} \quad
	{\rm corr}\{V(-c),V(kc)\}=-{3k\over (k+2)(2k+1)}. \eqno(4.4)$$ 
For small $h$ it is worth noting that 
$$\eqalign{
E\{V(b+h)-V(b)\}&={1\over b+h}-{1\over b}\doteq -{1\over b^2}h, \cr
E\{V(b+h)-V(b)\}^2&={2\over b^2}+{2\over (b+h)^2}-{4\over b(b+2h)}
		\doteq{4\over b^3}h. \cr}$$

\smallskip
{\sl 4C. Bivariate and multivariate exponential distributions.} 
We have constructed a full exponential process, 
and in particular $(V(b_1),\ldots,V(b_n))$ is a random vector
with dependent and exponential marginals. 
These bivariate and multivariate exponential classes of distributions
appear to be new. See Block (1985), for example, for a review
of the field of multivariate exponential distributions,
and see 6E below for a couple of other processes with 
exponential marginals.

Formula (4.4) shows that if values 
$\mu_1>0$, $\mu_2>0$, $\rho\in(0,1)$ are given, 
then a pair of dependent exponentials $(V(b_1),V(b_2))$ 
can be found with $EV(b_1)=\mu_1$, $EV(b_2)=\mu_2$, 
and correlation $\rho$. 
The class of bivariate exponential distributions is accordingly 
rich in the sense of achieving all positive correlations. 
The negative correlation in (4.4) 
starts out at zero for $k$ small, decreases to $-{1\over3}$ for $k=1$,
and then climbs up towards zero again when $k$ grows,
so negative correlations between $-{1\over3}$ and $-1$  
cannot be attained. Note that the maximal negative correlation 
occurs between $V(b)$ and $V(-b)$. 
%

In order to study the bivariate distribution for $(V(b),V(c))$ 
we calculate its double moment sequence (4.1) explicitly, 
for the case of $0<b<c$. 
The technique is to split the integral into $n!=(p+q)!$ parts, 
corresponding to all different orderings of the 
$n=p+q$ time indexes, the point being that a formula like (3.3) 
for the density of $(W(t_1),\ldots,W(t_n))$ can be exploited 
for each given ordering.  
These orderings can be grouped into 
$n\choose p$ types of paths, say $(e_1t_1,\ldots,e_nt_n)$ 
where $t_1<\cdots<t_n$ and 
$e_j$ is equal to $b$ in exactly $p$ cases 
and equal to $c$ in exactly $q$ cases. 
There are $p!q!$ different paths 
for given locations for the $p$ $b$'s and $q$ $c$'s, 
so the full integral can be written $\sum p!q!\,g({\rm path})$, 
where the sum is over all $n\choose p$ classes of paths and 
$g({\rm path})$ is the contribution for a specific path 
of the appropriate type. It remains to calculate the 
$g$-terms of various types, i.e.~to evaluate 
$$\int\cdots\int_{0<t_1<\cdots <t_n}
	f_{t_1,\ldots,t_n}(e_1t_1,\ldots,e_nt_n)\,dt_1\cdots dt_n $$
for a path with $e_j$'s equal to $b$ or $c$. 
Stameniforous integrations, 
similar to but more strenuous than those used to prove (4.3), 
show in the end that 
$$g({\rm path})
	=\Bigl({1\over b_0}\Bigr)^{i(0)}\Bigl({1\over b_1}\Bigr)^{i(1)} 
		\cdots \Bigl({1\over b_n}\Bigr)^{i(n)}, \eqno(4.5)$$
where the path when read backwards, i.e.~looking through
$(e_n,\ldots,e_1)$ in the notation above, has $i(0)$ $b$'s first, 
then $i(1)$ $c$'s, then $i(2)$ $b$'s, {\it\&}cetera. 
Furthermore $b_0=b$, $b_1=c$, and $b_j=b+j(c-b)$. 
Note that $i(0)+i(2)+\cdots=p$, $i(1)+i(3)+\cdots=q$, 
and $i(0)+i(1)+\cdots+i(n)=n$. And $EV(b)^pV(c)^q$ is equal to 
$p!q!$ times the sum of all such $g({\rm path})$ terms. 

To illustrate this somewhat cryptic formula, try $EV(b)^2V(c)^2$. 
There are ${4\choose2}=6$ types of paths, corresponding to 
$(b,b,c,c)$, 
$(b,c,b,c)$,
$(b,c,c,b)$,
$(c,b,b,c)$,
$(c,b,c,b)$,
$(c,c,b,b)$,
and each of these has weight $2!2!=4$. 
Their $(i(0),i(1),\ldots,i(4))$ representations are 
respectively 
$(0,2,2,0,0)$,
$(0,1,1,1,1)$,
$(1,2,1,0,0)$,
$(0,1,2,1,0)$,
$(1,1,1,1,0)$, 
$(2,2,0,0,0)$. 
Accordingly 
$$EV(b)^2V(c)^2=
4\Big\{{1\over b_1^2b_2^2}
+{1\over b_1b_2b_3b_4}
+{1\over b_0b_1^2b_2}
+{1\over b_1b_2^2b_3}
+{1\over b_0b_1b_2b_3}
+{1\over b_0^2b_1^2}\Bigr\},$$
where $b_0=b$, $b_1=c,\ldots,b_4=b+4(c-b)$. 

We have not been able to produce an explicit formula for the
joint probability density of $(V(b),V(c))$, 
but at least an expression can be found for 
its joint moment generating function. 
It becomes
$$\eqalign{ 
E\exp\{sV(b)+tV(c)\} 
&=\sum_{n\ge0}{1\over n!}\sum_{p+q=n}{n\choose p}
	s^pt^q\,E\{V(b)^pV(c)^q\} \cr
&=\sum_{n\ge0}{1\over n!}\sum_{p+q=n}{n!\over p!q!}s^pt^q\,
	\sum_{\rm paths}p!q!\,g({\rm path}) \cr
&=\sum_{p\ge0,\,q\ge0}s^pt^q\,
	\sum_{i(0),i(1),\ldots,i(p+q)}\Bigl({1\over b_0}\Bigr)^{i(0)} 
  	\Bigl({1\over b_1}\Bigr)^{i(1)} 
       	\cdots \Bigl({1\over b_{p+q}}\Bigr)^{i(p+q)}, \cr} \eqno(4.6)$$
where again the inner sum is over all $(p+q)!/p!q!$ types of paths
with $p$ $b$'s and $q$ $c$'s, and the multiplicities
$i(0),i(1),\ldots,i(p+q)$ have even-sum $p$ and odd-sum $q$, 
as explained above. 

One can similarly establish formulae for product moments 
of more than two $V(b)$'s, and investigate other aspects
of the multivariate exponential distributions associated with
the $V(.)$ process. We remark that these distributions can
be simulated, with some effort, through using $V_\eps(b)$ with
a small $\eps$, and this is one way of computing bivariate and
multivariate probabilities when needed. 
Another way would be via numerical inversion of the joint
moment generating function. 

\bigskip
{\bf 5. Total relative time along a line after time c.} 
As a generalisation of (1.2), consider the total relative time
spent along the ray $w=bt$ during $t\ge c$, i.e.
$$V_\eps(c,b)={1\over\eps}\int_c^\infty 
	I\{bt-\half\eps\le W(t)\le bt+\half\eps\}\,dt. \eqno(5.1)$$
The story told in the final paragraph of Section 1 is one motivation
for studying these variables.   
The main result about them is that 
$$V_\eps(c,b)\rightarrow_d V(c,b)\sim k(|b|\sqrt{c})\,{\rm Exp}(|b|)
		+(1-k(|b|\sqrt{c}))\,\delta_0, \eqno(5.2)$$
where again $\delta_0$ is degenerate at zero 
and $k(u)=2(1-\Phi(u))$.
Note that $k(|b|\sqrt{c})=1$ when $c=0$, 
so that (5.2) indeed contains our earlier result (3.1). 

It is possible to prove this by establishing 
a differential equation for the Laplace transform
of $V(c,b)$ with appropriate boundary conditions, 
and then solve, as in Section 2,
but it is as convenient to prove moment convergence.
Take $b>0$ for simplicity. It takes one moment to show that 
$$EV_\eps(c,b)
	\rightarrow \int_c^\infty f_t(bt)\,dt  
	=\int_c^\infty \phi(b\sqrt{t})/\sqrt{t}\,dt
	={1\over b}\int_{b\sqrt{c}}^\infty 2\phi(x)\,dx
	=k(b\sqrt{c})/b. $$
And when $p\ge2$ we find 
$$\eqalign{
EV_\eps(c,b)^p&\rightarrow
	p!\int\cdots\int_{c\le t_1<\cdots<t_p}
	f_{t_1,\ldots,t_p}(bt_1,\ldots,bt_p)\,dt_1\cdots dt_p \cr
	&=p!\int_c^\infty \int_0^\infty\cdots\int_0^\infty
	 {\phi(b\sqrt{u_1})\over \sqrt{u_1}}\cdots 
	 {\phi(b\sqrt{u_p})\over \sqrt{u_p}}\,du_1\cdots du_p \cr
	&=p!\,{k(b\sqrt{c})\over b}
		\Bigl({1\over b}\Bigr)^{p-1}. \cr} \eqno(5.3)$$
The Laplace transform function of 
this limit distribution candidate becomes 
$$\eqalign{
E\exp(-\lambda V)&=1+\sum_{p=1}^\infty
	{(-\lambda)^p\over p!}p!{k(b\sqrt{c})\over b}
	\Bigl({1\over b}\Bigr)^{p-1} \cr
	&=1+k(b\sqrt{c}){-\lambda/b \over 1+\lambda/b} 
	 =k(b\sqrt{c}){b\over b+\lambda}+1-k(b\sqrt{c}), \cr}$$
which is recognised as the moment generating function
of the mixture variable that with probability $k(b\sqrt{c})$ is
an exponential with parameter $b$ 
and with probability $1-k(b\sqrt{c})$ is equal to zero.
This proves (5.2). \firkant

\smallskip
{\csc Remark.} 
Let us briefly discuss a specific consequence, namely that 
${\rm Pr}\{V_\eps(c,b)=0\}$ in this situation converges to 
${\rm Pr}\{V(c,b)=0\}$, which is 
$1-k(b\sqrt{c})=2\Phi(b\sqrt{c})-1$. 
But having $V_\eps(c,b)=0$ in the limit 
means that $W(t)$ stays away from $bt$ during $[c,\infty)$, 
and it cannot stay above the curve all the time
since $W(t)/t$ goes to zero. Hence $2\Phi(b\sqrt{c})-1$ is simply 
the probability that $W(t)<bt$ during all of $[c,\infty)$,
or ${\rm Pr}\{\max_{t\ge c}W(t)/t<b\}$. 
Using finally the transformation
$W^*(t)=tW(1/t)$ to another Brownian motion one sees that 
$${\rm Pr}\{\max_{0\le t\le 1/c}W(t)\le b\}
	=2\Phi(b\sqrt{c})-1={\rm Pr}\{|W(1/c)|\le b\}. \eqno(5.4)$$
We have in other words rederived a classic 
distributional result for Brownian motion.~\firkant

\smallskip
The distribution of $V(c,-1)-V(c,1)$ 
comes up in the statistical estimation problem discussed
in Section 1; see also 6C below and Hjort and Fenstad (1992b, Section 6). 
When $c=0$ this is a difference
between two unit exponentials with intercorrelation $-{1\over3}$. 
The case $c>0$ is more complicated. Then 
$$\bigl(V(c,-1),V(c,1)\bigr)=\cases
	{0 &with probability $\pi_{00}$, \cr
	(U_{-1},0) &with probability $\pi_{10}$, \cr
	(0,U_1)  &with probability $\pi_{01}$, \cr
	(U_{-1},U_1) &with probability $\pi_{11}$, \cr} \eqno(5.5)$$
in which $U_{-1}$ and $U_1$ are unit exponentials with a certain
dependence structure. 
Furthermore $\pi_{00}$ is the probability that 
$W(t)$ stays between $-t$ and $t$ during $[c,\infty)$,
$\pi_{10}$ is the probability that $W(t)$ comes below $-t$ but is 
never above $t$, 
$\pi_{01}$ is the probability that $W(t)$ comes above $t$ but is 
never below $-t$, and 
$\pi_{11}$ is the probability that $W(t)$ experiences both $W(t)<-t$ 
and $W(t)>t$ during $[c,\infty)$. 
When $c=0$ then $\pi_{11}$ is 1 and the others are zero. 
In the positive case these probabilities can be found in terms of 
$H(u)$, the probability that $\max_{0\le s\le 1}|W(s)|\le u$, 
by the transformation arguments used to reach (5.4). 
One finds 
$$\pi_{00}=H(\sqrt{c}), \quad
	\pi_{01}=\pi_{10}=2\Phi(\sqrt{c})-1-\pi_{00}, \quad
	\pi_{11}=1-\pi_{00}-\pi_{01}-\pi_{10}, $$
in which $H(u)={\rm Pr}\{\max_{0\le s\le 1}|W(s)|\le u\}$.
A classic alternating series expression for $H(u)$ 
can be found in Shorack and Wellner 
(1986, Chapter 2, Section 2), for example,
and a new way of deriving this formula is by calculating 
all product moments $EV(-c)^pV(c)^q$ and then study the analogue of 
(4.6). This would be analogous to the way in which (5.4) was proved 
above, but the present case is much more laborious. 
Here we merely note that 
$$EV(c,-1)V(c,1)=\pi_{11}EU_{-1}U_1=\hbox{$2\over3$}k(3\sqrt{c}),$$
from which the correlation between $U(-1)$ and $U(1)$ also can be read off. 

\bigskip
{\bf 6. Supplementing results.} 

\smallskip
{\sl 6A. Total relative time along the time axis.}
The variable $V_\eps(b)$ of (1.2) is infinite when $b=0$. 
But consider
$$V_{\eps,T}={1\over\eps}{1\over\sqrt{T}}
	\int_0^T\{-\half\eps\le W(t)\le \half\eps\}\,dt, \eqno(6.1)$$
the relative time along the time axis during $[0,T]$. 
The moment sequence converges as $\eps\rightarrow0$ and $T\rightarrow\infty$,
as follows, using (3.3) once more: 
$$\eqalign{
E(V_{\eps,T})^p
	&={p!\over \eps^p\,T^{p/2}}
	  \int\cdots\int_{0<t_1<\cdots<t_p<T}
	  \Bigl[\phi(0)^p{\eps\over \sqrt{t_1}}\cdots
		{\eps\over \sqrt{t_p-t_{p-1}}}
	 	+O(\eps^{p+1})\Bigr]\,dt_1\cdots dt_p \cr
	&\rightarrow p!\,\phi(0)^p\int\cdots\int_{0<x_1<\cdots<x_p<1}
	  x_1^{-1/2}\cdots (x_p-x_{p-1})^{-1/2}(1-x_p)^0\,dx_1\cdots dx_p \cr
	&={p!\over (2\pi)^{p/2}}{\Gamma(\half)^p\Gamma(1)\over \Gamma(p/2+1)}
	 ={1\over 2^{p/2}}{p!\over \Gamma(p/2+1)}. \cr}$$
The limit distribution candidate $V_0$ has consequently 
$EV_0^{2p}=(\half)^p(2p)!/p!$, which means that $V_0^2$ gets 
moment generating function $(1-2t)^{-1/2}$. 
So $V_0^2$ is a $\chi^2_1$ (since the distribution of a 
chi-squared is determined by its moments), i.e.~$V_0$ is a $|N(0,1)|$. 
	  
It wasn't necessary here to send $T$ to infinity,
since the scaling property for $W$ 
[$W^*(t)=W(ct)/\sqrt{c}$ gives a new Brownian motion]
implies that the limit distribution of $V_{\eps,T}$ 
as $\eps\rightarrow0$ is independent of $T$.

One generalisation of this is in the following direction. 
Instead of (6.1), look at 
$$V_{\eps,T}={1\over \eps}{1\over \sqrt{T}}\int_0^T
	h\bigl(W(t)/\eps)\bigr)\,dt
	={1\over \sqrt{T/\eps^2}}\int_0^{T/\eps^2} h(W^*(t))\,dt, $$
where $h(x)$ is any function with bounded support, and where 
$W^*$ in the second expression is another Brownian motion 
obtained from the first one by transformation. 
The case considered earlier is $h(x)=I\{|x|\le\half\}$. 
It can be shown that $V_{\eps,T}\rightarrow a|N(0,1)|$ 
in distribution as $T/\eps^2\rightarrow\infty$, 
where $a$ is a constant depending on $h$.
This is not easy to prove via the moment convergence
technique, but can be established using methods from 
Khasminskii (1980). 
 

\smallskip
{\sl 6B. Implications for partial sum processes.}
Let us first point out that an alternative construction of 
our total relative time variables is to use $I\{bt\le W(t)\le bt+\eps\}$
instead of $I\{W(t)\in bt\pm\half\eps\}$ in (1.2) and (5.1). 
Results of previous sections hold equally for this alternative
definition of $V_\eps(b)$ and $V_\eps(c,b)$, and this 
is a bit more convenient in 6C below. 
Now suppose $X_1,X_2,\ldots$ are i.i.d.~with 
mean $\xi$ and variance $\sigma^2$,
and consider the normalised partial sum process 
$W_m(t)=m^{-1/2}\sum_{i=1}^{[mt]}(X_i-\xi)/\sigma$. 
In particular $W_m({n\over m})=S_n/\sqrt{m}$,
writing $Y_i=(X_i-\xi)/\sigma$ and $S_n$ for their partial sums, 
and $W_m(.)$ converges to Brownian motion by Donsker's theorem.
Motivated by (1.2) and (5.1) we define
$$\eqalign{
V_{m,\eps}(c,b)&={1\over\eps}{1\over m}\sum_{n/m\ge c}
	I\Bigl\{b{n\over m}\le {S_n\over \sqrt{m}}
			\le b{n\over m}+\eps\Bigr\} \cr
&={1\over\eps}{1\over m}\sum_{n/m\ge c}
	I\Bigl\{b\sqrt{n\over m}\le T_n=\sqrt{n}(\bar X_n-\xi)/\sigma
		\le b\sqrt{n\over m}+\eps\sqrt{m\over n}\Bigr\} \cr
&={1\over\eps}\int_{\langle cm\rangle/m}^\infty
       	I\Bigl\{b{[mt]\over m}\le W_m(t)
		\le b{[mt]\over m}+\eps\Bigr\}\,dt, \cr} \eqno(6.2)$$
where $\langle cm\rangle$ 
denotes the smallest integer exceeding or equal to $cm$. 
It is clear that this variable is close to $V_\eps(c,b)$ for large $m$,
and should accordingly converge in distribution to $V(c,b)$ of (4.2) 
when $m\rightarrow\infty$ and $\eps\rightarrow0$. 

\smallskip
{\csc Proposition.} 
{\sl Assume that the $X_i$'s have a finite third absolute moment. 
If $c>0$ is fixed, then $V_{m,\eps(m)}(c,b)\rightarrow_d V(c,b)$
if only $\eps(m)\rightarrow0$ as $m\rightarrow0$. And 
$$V_{m,\eps(m)}(c(m),b)\rightarrow_d V(b)\sim {\rm Exp}(|b|) \eqno(6.3)$$ 
provided 
$\eps(m)\rightarrow0$,
$c(m)\rightarrow0$, 
$mc(m)\rightarrow\infty$, and 
$\eps(m)/c(m)^{1/2}\rightarrow0$.
\smallskip}

{\csc Proof:}
This can be proved in various ways and under various conditions.
One feasible possibility is to demonstrate moment convergence 
of $E\{V_{m,\eps(m)}(c,b)\}^p$ towards the right hand side of (5.3),
for each $p$. 
One basically needs the smallest $n$ 
in the sum to grow towards infinity, 
so that the central limit theorem and Edgeworth--Cram\'er expansions 
can begin to work, 
and the largest of all $\eps(m)\sqrt{m/n}$ terms to go to zero, 
so that Taylor expansions can begin to work;
see the middle term in (6.2). When $c$ is fixed then the sum is 
over all $n\ge mc$, and it suffices to have $\eps(m)\rightarrow0$ as 
$m\rightarrow\infty$. To reach $V(b)=V(0,b)$ in the limit we 
need the stated behaviour for $\eps(m)$ and $c(m)$. 
We have used the third moment assumption to 
bound the error $r(t)$ in the Edgeworth expression 
$G_n(t)=\Phi(t)+r(t)$ for the distribution of 
$T_n=\sqrt{n}(\bar X_n-\xi)/\sigma$; 
one has $|r_n(t)|\le cn^{-1/2}/(1+|t|)^3$, 
and this is helpful when it comes to verifying  
conditions when employing 
Lebesgue's theorem on dominated convergence.~\firkant


\smallskip
We may conclude that the total relative time along $b{n\over m}$ 
for the normalised partial sum process has a limit distribution, 
which is either exponential or of the mixture type (4.2). 
The middle expression also invites $V_{m,\eps}$ to be thought of
as the total relative time for the normalised $T_n$ process 
along the square root boundary $b\sqrt{n/m}$. 
The result is also valid for 
$T_n=\sqrt{n}(\theta_n-\theta)/\sigma$ 
in a more general estimation theory setup; 
see Hjort and Fenstad (1992a, 1992b).   

The 6A result has also implications for partial sum processes.
One can prove that 
$${1\over \eps}\int_0^1I\{|W_m(t)|\le \half\eps\}\,dt
	\rightarrow_d |N(0,1)|$$
when $\eps\rightarrow0$ and $m\rightarrow\infty$, under suitable
conditions. This implies for example that 
$m^{-1/2}\sum_{i=1}^mI\{|S_i|\le\half\}$ has the absolute normal limit,
as does $m^{-1/2}\sum_{i=1}^mI\{S_i=0\}$ 
for the random walk process. 

\smallskip
{\sl 6C. Second order asymptotics for the number of $\delta$-errors.}
To show how the total relative time variables for Brownian motion
are related to the estimation theory problems described in Section 1, 
consider the structurally simple case of i.i.d.~variables
$X_i$ with mean $\xi$ and standard deviation $\sigma$, 
and where ${n\over n+k}\bar X_n$ 
is used to estimate $\xi$. Consider $Q_\delta(k)$, the number of times
$|{n\over n+k}\bar X_n-\xi|\ge\delta$,
counted among $n\ge c/\delta^2$. Then $\delta^2Q_\delta(k)$ tends
to $Q=Q(c,1/\sigma)$ of (1.2), for each choice of $k$,
and $\delta^2$ times $Q_\delta(k)-Q_\delta(0)$ goes to zero. 
This follows from results in Hjort and Fenstad (1992a). 
But $\delta\{Q_\delta(k)-Q_\delta(0)\}$ can be written $A_\delta-B_\delta$,
after some analysis, where 
$$\eqalign{
A_\delta&=\sqrt{m}\int_{\langle mc\rangle/m}^\infty
	I\Bigl\{-{[mt]\over m}{1\over \sigma}\le W_m(t) \le
	-{[mt]\over m}{1\over\sigma}+{1\over \sqrt{m}}{k\xi\over \sigma}
		-{1\over m}{k\over \sigma}\Bigr\}\,dt, \cr
B_\delta&=\sqrt{m}\int_{\langle mc\rangle/m}^\infty
	I\Bigl\{{[mt]\over m}{1\over \sigma}\le W_m(t) \le
	{[mt]\over m}{1\over\sigma}+{1\over \sqrt{m}}{k\xi\over \sigma}
		+{1\over m}{k\over \sigma}\Bigr\}\,dt, \cr}$$
and where $m=1/\delta^2$. 
These variables resemble those considered in (6.2) and (6.3). 
With $\eps=m^{-1/2}k\xi/\sigma$ we have 
$$A_\delta\doteq_d k\xi/\sigma\,V_{m,\eps}(c,-1/\sigma)
	\quad {\rm and} \quad
  B_\delta\doteq_d k\xi/\sigma\,V_{m,\eps}(c,1/\sigma),$$
where `$\doteq_d$' signifies that the difference goes to zero in probability.
It follows from the result of 6B that 
$$\delta\{Q_\delta(k)-Q_\delta(0)\}\rightarrow_d
	k\xi/\sigma\,\{V(c,-1/\sigma)-V(c,1/\sigma)\}
	\quad {\rm as\ }\delta\rightarrow0. \eqno(6.4)$$
This is also true with $c=0$ in the limit, i.e.~with 
$k\xi/\sigma\,\{V(-1/\sigma)-V(1/\sigma)\}$ on the right hand side, 
provided $c=c(\delta)=\delta$ is used in the definition of 
$Q_\delta(k)$ and $Q_\delta(0)$. 
Note the relevance of (5.5) for the present problem. 

Hjort and Fenstad (1992b) also work with the direct 
expected value of $Q_\delta(k)-Q_\delta(0)$ and similar variables.
These converge to explicit functions of $k$ 
(and other parameters, in more general situations), 
which can then be minimised to single out estimator sequences with the 
second order optimality property of having the smallest expected 
number of $\delta$-errors. 
This is done in Hjort and Fenstad (1992b), 
in several situations. 
We remark that the skewness 
$\gamma=E(X_i-\xi)^3/\sigma^3$ is not involved in 
(6.4), but is prominently present in the limit of 
$E\{Q_\delta(k)-Q_\delta(0)\}$, and its minimisation. 

\smallskip
{\sl 6D. Relative time along other curves.} 
To generalise our framework, consider 
$$V_\eps={1\over \eps}\int_a^\infty
	I\{b(t)\le X(t)\le b(t)+\eps g(t)\}\,dt, \eqno(6.5)$$
where $x=b(s)$ is some curve of interest 
and $g(t)$ a possible scaling factor. In many cases there is a 
distributional limit as $\eps\rightarrow0$, and perhaps 
the first couple of moments can be obtained. 
The limit distribution is simple only for cases that can be 
transformed back to (1.1) and (5.1), however. 
For an example, we note that 
the total relative time an Ornstein--Uhlenbeck process 
$X(t)$ spends along $be^t$ can be shown to be exponential,
for example, with suitable $g(t)$ in (6.5). 

\smallskip
{\sl 6E. Other exponential and gamma processes.} 
(i) 
If $U(b)=|b|V(b)$, then $U(b)$ is unit exponential 
for each $b$. In particular its marginal 
mean and variance are constant, 
and ${\rm cov}\{U(b),U(b+h)\}=b/(b+2h)$. 
(ii) 
By adding independent copies of $V(.)$ (or $U(.)$)
we get processes with marginals that are gamma distributed. 
This leads in particular to bivariate and 
multivariate gamma distributions or chi-squared distributions
(with even-numbered degrees of freedom only). 
(iii) 
There are other processes that share with $V$ and $U$ 
the property of having exponentially distributed marginals. 
One example is $V^*(b)=\half\{W_1(b)^2+W_2(b)^2\}$, 
where $W_1$ and $W_2$ are independent Brownian motions. 
This is a Markov process, while our $V(b)$ process is not. 
The possible correlations of $(V^*(b_1),\ldots,V^*(b_n))$
span a smaller space than those of $(V(b_1),\ldots,V(b_n))$,
indicating that the $V^*$ process may be less adequate when it comes
to building multivariate exponential models. 
(iv) 
And yet another process with exponential marginals is provided by 
$V^{**}(b)=\max_{t\ge0}\{W(t)-bt\}$, defined for positive $b$.
This process is studied by \c Cinlar (1992). 
If $X(t)=W(t)-bt$, then \c Cinlar's paper is concerned with the maximum
value of this process and where the maximum occurs, whereas
the present paper has been concerned with the amount of time
such a process spends along a line. 

\bigskip
{\bf Acknowledgments.}
This paper was written 
at the Mathematical Sciences Research Institute at Berkeley,
where we spent some pleasant weeks as invited participants in
the 1991 Statistics Programs.  
This was made possible through 
generous support from National Science Foundation 
grant 8505550. 

\bigskip
\parindent0pt\parskip3pt\baselineskip12pt
\centerline{\bf References}

\smallskip
%
%

Block, H.W., 
Multivariate exponential distributions,
{\sl Encyclopedia of Statistical Science,}
Eds.~Kotz and Johnson
(Wiley, New York), 1985. 

%

\c Cinlar, E., 
Sunset over Brownistan,
{\sl Stochastic Processes and their Applications} {\bf 40}, 
1992, 45--53. 

Hjort, N.L.~and Fenstad, G., 
On the last time and the number of times an estimator
is more than $\eps$ from its target value,
{\sl Annals of Statistics} {\bf 20}, 1992a, 469--489. 

Hjort, N.L.~and Fenstad, G., 
Some second order asymptotics for the number of times an 
estimator is more than $\eps$ from its target value,
Statistical research report, University of Oslo, 1992b.

It\^o, K.~and McKean, H.P.~Jr., 
{\sl Diffusion Processes and their Sample Paths} 
(2nd edition, Springer-Verlag, Berlin), 1979.

Karlin, S.~and Taylor, H.M., 
{\sl A Second Course in Stochastic Processes}
(Academic Press, Toronto), 1981. 

Khasminskii, R.Z., 
{\sl Stability of Stochastic Differential Equations}
(Sijthoff and Noordhoff, Amsterdam), 1980.


Shorack, G.R.~and Wellner, J.A., 
{\sl Empirical Processes with Applications to Statistics}
(Wiley, New York), 1986. 

\bye